\documentclass{article}

\usepackage[T1]{fontenc}
\usepackage[english]{babel}
\usepackage{amsmath}
\usepackage{amsthm}
\usepackage{amssymb}

\newcommand{\comments}[1]{} 
\newcommand{\leftsup}[2]{{}^{#1}\!{#2}}	

\newtheorem{thm}{Theorem}

\newtheorem{lem}{Lemma}
\newtheorem{cor}{Corollary}

\newtheorem{prop}{Proposition}

\theoremstyle{definition}

\theoremstyle{remark}
\newtheorem{rem}{Remark}

\theoremstyle{remark}
\newtheorem{ex}{Example}

\def\cI{\mathcal I}

\def\cM{\mathcal M}

\def\cP{\mathcal P}
\def\cC{\mathcal C}

\def\cH{\mathcal H}
\def\cV{\mathcal V}

\def\cB{\mathcal B}

\def\bR{\mathbb R} 

\def\bN{\mathbb N}

\def\bC{\mathbb C}

\def\xx{\mathbf x}
\def\yy{\mathbf y}

\def\Harm{\mathit{Harm}}

\def\Ker{\operatorname{Ker}}

\def\mand{\text{\ \ \ and\ \ \ }}

\def\clifford{\cC\ell}

\def\pa{\partial}

\def\bul{\;\bullet}
\def\wed{\;\wedge}
\def\pod{\underline}
\def\nad{\overline}

\def\eX{\leftsup{e}{X}}
\def\oX{\leftsup{o}{X}}

\def\pyy{\underline{\yy}}

\begin{document}

\title{The Gelfand-Tsetlin bases for Hodge-de Rham systems in Euclidean spaces}

\author{Richard Delanghe, Roman L\' avi\v cka and Vladim\'ir Sou\v cek}

\maketitle

\begin{abstract}
The main aim of this paper is to construct explicitly orthogonal bases for the spaces $\cH^s_k(\bR^m)$ of $k$-homogeneous polynomial solutions of the Hodge-de Rham system in the Euclidean space $\bR^m$ which take values in the space of $s$-vectors. Actually, we describe even the so-called Gelfand-Tsetlin bases for such spaces in terms of Gegenbauer polynomials.
As an application, we obtain an algorithm how to compute an orthogonal basis of the space of homogeneous solutions of a generalized Moisil-Th\'eodoresco system in $\bR^m.$  

\medskip\noindent
{\bf Keywords:} Clifford analysis, Hodge-de Rham  system, Gelfand-Tsetlin basis, generalized Moisil-Th\'eodoresco system

\medskip\noindent
{\bf AMS classification:} 30G35, 58A10, 22E70, 33C50 
\end{abstract}


\section{Introduction}
   
In what follows, we are interested mainly in the spaces $\cH^s_k(\bR^m)$ of $k$-homogeneous polynomial solutions of the Hodge-de Rham system in the Euclidean space $\bR^m$ which take values in the space of $s$-vectors.
As is well known (see \cite{hom, DLS}), the space $\cH^s_k(\bR^m)$ can be viewed naturally as an irreducible finite dimensional module over the orthogonal group $O(m).$ 
Moreover, even for any irreducible finite dimensional module even over a general classical simple Lie algebra, an abstract definition of the Gelfand-Tsetlin (GT for short) basis is given, see e.g.\ \cite{mol, GT}. 
The main aim of this paper is to describe an explicit construction of GT bases for the spaces $\cH^s_k(\bR^m).$ 
Let us emphasize that the GT basis is always orthogonal with respect to any invariant inner product on the given module. 
For the spaces $\cH^s_k(\bR^m),$ 
bases (not always orthogonal) were so far constructed and studied only in some special cases (e.g., for $s=1$ or in low dimensions),
see \cite{cno, R1, R2, cac, leu, zei, mor09,GM10b}.
We work within the frame of Clifford analysis but all the results of the paper can be easily translated into the language of differential forms, as is explained in \cite{BDS}.

Clifford analysis can be considered as a~refinement of harmonic analysis (see e.g.\ \cite{BDS2, DSS, GM}). It studies mainly solutions of the Dirac equation in $\bR^m$ which take values in the Clifford algebra $\clifford_m$ over $\bR^m.$
As such solutions are real analytic an important step is to understand first the structure of homogeneous polynomial solutions.
On the space of $\clifford_m$-valued polynomials in $\bR^m,$ we can consider the so-called $L$-action.
The role of building blocks are then played by the spaces of homogeneous spinor valued solutions which are irreducible in this case. An explicit construction of orthogonal (or even GT) bases for these spaces in any dimension was explained in \cite{BGLS}. 
Let us remark that, for the construction, the so-called Cauchy-Kovalevskaya (CK for short) method developed already in \cite[Theorem 2.2.3, p. 315]{DSS} was used. 
In the introduction to the paper \cite{BGLS},
more details on history
of this topic are available.
For the classical case we can refer further to \cite{Bock2010c,Bock2010a,Bock2009,BG,CacGueMal,CacGueBock,BockCacGue,cac,CM06,CM07,CM08,
FCM,FM,Gurlebeck1999,NGue2009,lavSL2,step2,Malonek1987,MS,mor09,som,van}. 
Analogous results in Hermitean Clifford analysis are described in \cite{ckH, kerH,GTinH, GT2H}. 

On the other hand, on the space of Clifford algebra valued polynomials the $H$-action, defined in \eqref{Haction} below,
can be considered. 
In a series of papers \cite{R1,R2,DLS,DLS3,lav_isaac09}, generalized Moisil-Th\'eodoresco (GMT for short) systems have been studied. The spaces of homogeneous solutions of a GMT system are important examples of modules under the $H$-action. More explicitly, 
let $S\subset\{0,1,\ldots,m\}$ be given and put $$\clifford_m^S=\bigoplus_{s\in S}\clifford_m^s,$$ 
where $\clifford_m^s$ is the space of $s$-vectors in $\clifford_m.$ 
Moreover, denote by $\cM^S_k(\bR^m)$
the set of $k$-homogeneous solutions of the Dirac equation in $\bR^m$ (that is, spherical monogenics) which take values in $\clifford_m^S.$
It is clear that, for $S=\{0,1,\ldots,m\},$ the space $\cM^S_k(\bR^m)$ coincides with the space 
of all $k$-homogeneous Clifford algebra valued spherical monogenics in $\bR^m.$ Moreover, for $S=\{s\},$  $\cM^S_k(\bR^m)=\cH^s_k(\bR^m).$
Actually, it turns out that basic building blocks for the space $\cM^S_k(\bR^m)$ are just the spaces $\cH^s_k(\bR^m)$ 
of homogeneous solutions of the Hodge-de Rham system (see Theorem \ref{tGMT}). 

The main result of this paper is an explicit construction of orthogonal (or even GT) bases of the spaces $\cH^s_k(\bR^m).$ In Section 2, we recall briefly the Fischer decomposition for the $H$-action 
and the notion of GT basis at least for the orthogonal groups. 
We can construct the GT bases of the spaces $\cH^s_k(\bR^m)$ using  
the CK method analogously as in the case of spinor valued spherical monogenics. 
But we need first to adapt the CK method for this case, which is done in Section 3.
Finally, in Section 4, we describe an algorithm how to express GT bases 
of the spaces $\cH^s_k(\bR^m)$ in terms of Gegenbauer polynomials by induction on the dimension $m$ (see Theorem \ref{induction}). Moreover, we give explicit examples of these bases in dimensions 3 and 4 at least for small values of $k.$
As an application, we obtain an algorithm how to compute explicitly an orthogonal basis of the space $\cM^S_k(\bR^m)$ of homogeneous solutions for an arbitrary GMT system. In this paper, we just describe a~construction of these bases. A~detailed study of properties of the constructed bases will be given in a~next paper.

\section{Notations and known facts}

Let the vectors $e_1,\ldots,e_m$ form the standard basis of the Euclidean space $\bR^m.$
Denote by $\bR_{0,m}$ the real Clifford algebra over $\bR^m$ satisfying the basic multiplication relations
$$e_ie_j+e_je_i=-2\delta_{ij}$$
and by $\bC_m$ the corresponding complex Clifford algebra. 
In what follows, $\clifford_m$ means either $\bR_{0,m}$ or $\bC_m.$
 
The Clifford algebra $\clifford_m$ can be viewed naturally as the graded associative algebra
$$\clifford_m=\bigoplus_{s=0}^m\clifford_m^s.$$
Here $\clifford_m^s$ stands for the space of $s$-vectors in $\clifford_m.$
As usual, we identify a~vector $(x_1,\ldots,x_m)$ of $\bR^m$ with the 1-vector $x_1e_1+\cdots+x_m e_m$ of $\clifford_m^1.$ 
For a 1-vector $u$ and an $s$-vector $v,$ the Clifford product $uv$ splits into the sum of an $(s-1)$-vector $u\bullet v$ and an $(s+1)$-vector $u\wedge v.$ Indeed, we have that
$$uv=u\bullet v+u\wedge v\text{\ \ \ with\ \ \ }u\bullet v=\frac 12(uv-(-1)^svu)\text{\ \ and\ \ }u\wedge v=\frac 12(uv+(-1)^svu).$$
By linearity, we extend the so-called inner product $u\bullet v$ and the outer product $u\wedge v$ for a 1-vector $u$ and an arbitrary Clifford number $v\in\clifford_m.$

In what follows, we deal with the space $\cP^*=\cP^*(\bR^m)$ of $\clifford_m$-valued polynomials in the vector variable $\xx=(x_1,\ldots,x_m)$ of $\bR^m.$
Denote by $\cP^*_k$ the space of $k$-homogeneous polynomials of $\cP^*$ and by $\cP^s_k$ the space of $s$-vector valued polynomials of $\cP^*_k.$
In general, for $\cV\subset\cP^*$ put
$\cV_k=\cV\cap\cP^*_k$ and $\cV^s_k=\cV\cap\cP^s_k.$ 

\paragraph{The Fischer decomposition for the $H$-action}
On the space $\cP^*$ of Clifford algebra valued polynomials we can consider 
the so-called $H$-action of the Pin group $Pin(m),$ given by
\begin{equation}\label{Haction}
[H(r)P](\xx)=r\;P(r^{-1}\xx\; r)\;r^{-1},\ r\in Pin(m),\ P\in\cP^*\text{\ \ and\ \ }\xx\in\bR^m.
\end{equation}
Recall that the group $Pin(m)$ is a~double cover of the orthogonal group $O(m).$ 
Obviously, the multiplication by the vector variable $\xx=e_1x_1+\cdots+e_m x_m$ and the Dirac operator 
$$\pa =e_1\pa_{x_1}+\cdots+e_m\pa_{x_m}$$ (both applied from the left)
are examples of invariant linear operators on the space $\cP^*$ with the $H$-action.
On the other hand, we can split the left multiplication by a 1-vector $\xx$ into the outer multiplication $(\xx\wed)$ and the inner multiplication $(\xx\bul),$ that is, $$\xx=(\xx\wed)+(\xx\bul) .$$
Analogously,
the Dirac operator $\pa$ can be split also into two parts
$\pa=\pa^++\pa^-$ where
$$\pa^+P=\sum_{j=1}^m e_j\wedge(\pa_{x_j}P)\text{\ \ and\ \ }\pa^-P=\sum_{j=1}^m e_j\bullet(\pa_{x_j}P).$$
Actually,
the operators $\pa^+,$ $\pa^-,$ $(\xx\wed)$ and $(\xx\bul)$ are, in a~certain sense, basic invariant operators for the $H$-action  (see \cite{DLS3} for more details). 
Moreover, denote by $\cH^s_k=\cH^s_k(\bR^m)$ the space of polynomials $P\in\cP^s_k$ satisfying
the Hodge-de Rham system of equations
\begin{equation}\label{HdR}
\pa^+P=0,\ \pa^-P=0.
\end{equation}
It is easily seen that
$$\cH^s_k=\{P\in\cP^s_k|\ \pa P=0\}.$$
Let $\Omega$ be the set of all non-trivial words in the letters $(\xx\wed)$ and $(\xx\bul).$ 
Note that $(\pa^+)^2=0,$ $(\pa^-)^2=0,$ $(\xx\wed)^2=0$ and $(\xx\bul)^2=0.$ In particular, we have that the set $\Omega$ 
looks like
\begin{equation}\label{Omega}
\Omega=\{1,\ (\xx\wed),\ (\xx\bul),\ (\xx\wed)(\xx\bul),\ (\xx\bul)(\xx\wed),\ (\xx\wed)(\xx\bul)(\xx\wed),\ \ldots\}.
\end{equation}
Then the right analogue of the Fischer decomposition for the $H$-action reads as follows (see \cite{DLS2,lavSL3}).

\begin{thm}\label{tFischerH}
The space $\cP^*$ of $\clifford_m$-valued polynomials in $\bR^m$ decomposes as
\begin{equation}\label{FischerH}
\cP^*=\bigoplus_{s=0}^{m}\bigoplus_{k=0}^{\infty}\bigoplus_{w\in\Omega}w\cH^s_k.
\end{equation}
\end{thm}

\begin{rem}
(i) In addition, we have that $\cH^s_k=\{0\}$ just for $s\in\{0,m\}$ and $k\geq 1.$ 
In the case when $\clifford_m=\bR_{0,m}$ (resp.\ $\bC_m$), we have that
$\cH^0_0=\bR$ (resp.\ $\bC$) and $\cH^m_0=\bR\, e_M$ (resp.\ $\bC\, e_M$) with $e_M=e_1e_2\cdots e_m.$ 
Moreover,
under the $H$-action, all non-trivial modules $\cH^s_k$ are irreducible and mutually inequivalent.

\medskip\noindent
(ii) It is easy to see that 
$w\cH^s_k=\{0\}$ if either $s=0$ and the word $w$ begins with the letter $(\xx\bul)$ or
$s=m$ and the word $w$ begins with the letter $(\xx\wed).$ Otherwise, each module $w\cH^s_k$ is equivalent to the module $\cH^s_k.$
\end{rem}

\paragraph{Invariant inner products}

Recall that, on each (finite-dimensional) irreducible representation of $Pin(m)$ there is always an invariant inner product determined uniquely up to a~positive multiple. In what follows, we describe two well-known realizations of the invariant inner product on the spaces  $\cH^s_k(\bR^m),$
namely, the $L_2$-inner product and the Fischer inner product. First, for each $P,Q\in \cP^*_k(\bR^m),$ we define the $L_2$-inner product of $P$ and $Q$ as
\begin{equation}\label{L2product}
(P,Q)_1=\int_{S^{m-1}}[\nad{P(\xx)}Q(\xx)]_0\;d\,\Sigma(\xx)
\end{equation}
where $S^{m-1}$ is the unit sphere in $\bR^m$ and $d\,\Sigma$ is the elementary surface element on $S^{m-1}.$ Here, for each Clifford number $a\in\clifford_m,$ $\nad a$ stands for its Clifford conjugate and $[a]_0$ for its scalar part. 

Now we introduce the Fischer inner product. Each $P\in\cP^*_k(\bR^m)$ is of the form
$$P(\xx)=\sum_{|\alpha|=k}a_{\alpha}\xx^{\alpha}$$
where the sum is taken over all multi-indexes $\alpha=(\alpha_1,\ldots,\alpha_m)$ of $\bN^m_0$ with  $|\alpha|=\alpha_1+\cdots+\alpha_m=k,$ all coefficients $a_{\alpha}$ belong to $\clifford_m$ and $\xx^{\alpha}=x_1^{\alpha_1}\cdots x_m^{\alpha_m}.$ For $P,Q\in\cP^*_k(\bR^m),$ we define the Fischer inner product of $P$ and $Q$ as
\begin{equation}\label{Fischerproduct}
(P,Q)_2=\sum_{|\alpha|=k}\alpha!\;[\nad a_{\alpha}b_{\alpha}]_0
\end{equation}
where $\alpha!=\alpha_1!\cdots\alpha_m!,$ $P(\xx)=\sum a_{\alpha}\xx^{\alpha}$ and $Q(\xx)=\sum b_{\alpha}\xx^{\alpha}.$
It is easily seen that
$$(P,Q)_2=[(\nad P(\frac{\pa\ }{\pa \xx})Q)(0)]_0
\text{\ \ \ with\ \ \ }\nad P(\frac{\pa\ }{\pa
\xx})=\sum_{|\alpha|=k}\nad a_{\alpha}\frac{\pa^{|\alpha|}}{\pa
\xx^{\alpha}}.$$ Here $\pa^{|\alpha|}/\pa
\xx^{\alpha}=(\pa^{\alpha_1}/\pa
x^{\alpha_1}_1)\cdots(\pa^{\alpha_m}/\pa x^{\alpha_m}_m)$ as usual.

\paragraph{The Gelfand-Tsetlin bases}
In this paper, we are interested in a~construction of GT bases for the spaces $\cH^s_k(\bR^m).$ 
It is well-known that, under the $H$-action, the spaces $\cH^s_k(\bR^m)$ are examples of irreducible modules with the highest weights consisting entirely of integers, see \cite{DLS}. Hence the spaces $\cH^s_k(\bR^m)$ can be viewed as irreducible modules over the orthogonal group $O(m).$
Let us briefly recall how to construct a GT basis for the given space $\cH^s_k(\bR^m).$ 

The first step consists in reducing the symmetry to the group $O(m-1),$ realized as the subgroup of orthogonal transformations of $O(m)$ fixing the last vector $e_m.$ It turns out that, under the action of the group $O(m-1),$ the space $\cH^s_k(\bR^m)$ is reducible and  
decomposes into a~multiplicity free direct sum of irreducible
$O(m-1)$-submodules
\begin{equation}\label{branch}
\cH^s_k(\bR^m)=\bigoplus_{\mu_{m-1}}\cH(\mu_{m-1}).
\end{equation}
Since this irreducible decomposition is multiplicity free the decomposition is obviously orthogonal with respect to any invariant inner product given on the module $\cH^s_k(\bR^m).$ Moreover, as an irreducible $O(m-1)$-module, each piece $\cH(\mu_{m-1})$ is uniquely characterized by its label consisting of the highest weight for the corresponding $SO(m-1)$-module and a~number of $\{0,\pm 1\}$ (see \cite{DLS} for the labels of the spaces $\cH^s_k(\bR^m)$). Hence we could use the label of $\cH(\mu_{m-1})$ as its index $\mu_{m-1}.$
Let us remark that the decomposition \eqref{branch} is a~special case of the so-called branching law from representation theory.

Of course, we can further reduce the symmetry to the group $O(m-2),$ the subgroup of orthogonal tranformations of $O(m)$ fixing the last two vectors $e_{m-1}, e_m.$
Then we can again decompose each piece $\cH(\mu_{m-1})$ of the decomposition (\ref{branch}) into irreducible $O(m-2)$-submodules $\cH(\mu_{m-1},\mu_{m-2})$ and so on. 

Hence we end up with the decomposition of the given $O(m)$-module
$\cH^s_k(\bR^m)$ into irreducible $O(2)$-modules $\cH(\mu).$
Moreover, any such module $\cH(\mu)$ is uniquely determined by the sequence of labels
\begin{equation}
\label{pattern}
\mu=(\mu_{m-1},\ldots,\mu_2).
\end{equation}
To summarize, we decompose the given module $\cH^s_k(\bR^m)$ into
the direct sum of irreducible $O(2)$-modules
\begin{equation}\label{branch+}
\cH^s_k(\bR^m)=\bigoplus_{\mu}\cH(\mu).
\end{equation}
Moreover, with respect to any given invariant inner product on the module $\cH^s_k(\bR^m),$ the decomposition (\ref{branch+}) is obviously orthogonal. 
Now it is easy to obtain an orthogonal basis of the space $\cH^s_k(\bR^m).$
Indeed, each irreducible $O(2)$-module $\cH(\mu)$ is either one-dimensional or two-dimensional. In the latter case,
the space $\cH(\mu)$ decomposes further as
$$\cH(\mu)=\cH(\mu^+)\oplus\cH(\mu^-)$$
where $\cH(\mu^{\pm})$ are one-dimensional $SO(2)$-modules with the highest weights $\pm j$ for some natural number $j.$ Hence we even get an orthogonal decomposition of the given module $\cH^s_k(\bR^m)$ into one-dimensional $SO(2)$-modules $\cH(\tilde\mu),$ $\tilde\mu\in P.$
Now we construct easily a~basis of the space $\cH^s_k(\bR^m)$ by taking a~non-zero vector $e(\tilde\mu)$ from each one-dimensional piece $\cH(\tilde\mu).$ 
The obtained basis $$E=\{e(\tilde\mu):\tilde\mu\in P\}$$ is called a~GT basis of the module $\cH^s_k(\bR^m).$
It is easily seen that the vector $e(\tilde\mu)$ is uniquely determined by its index $\tilde\mu$ up to a~scalar multiple. Moreover, by construction, the GT basis $E$ is orthogonal 
with respect to any invariant inner product, including the $L_2$-inner product \eqref{L2product} and the Fischer inner product \eqref{Fischerproduct}.

\section{The Cauchy-Kovalevskaya method}

As was explained, to construct explicitly the GT basis of an $O(m)$-module $\cH^s_k(\bR^m)$ it is  first necessary to decompose the module $\cH^s_k(\bR^m)$ into irreducible $O(m-1)$-submodules, cf. \eqref{branch}. Now we show that such a~decomposition can be obtained using the Cauchy-Kovalevskaya method (CK for short). 

So let a~polynomial $p$ of $\cH^s_k(\bR^m)$ be given. Then $p$ is a $\clifford_m^s$-valued $k$-homogeneous polynomial 
in the vector variable $\xx\in\bR^m$ which solves the Dirac equation $\pa p=0.$ In what follows, we split the vector variable $\xx$ of $\bR^m$ 
into the first $m-1$ variables $\pod\xx=e_1x_1+\cdots+e_{m-1}x_{m-1}$ and the last one $x_m.$ 
Moreover, put $$\pod\pa=e_1\frac{\pa\ }{\pa
x_1}+\cdots+e_{m-1}\frac{\pa\ }{\pa x_{m-1}}.$$
As is well-known, the CK  extension operator $CK=e^{e_m x_m\pod\pa}$ reconstructs the monogenic polynomial $p(\xx)$ from the initial polynomial $p_0(\pod \xx)=p(\pod \xx,0),$ that is, $p=CK(p_0).$ In Theorem \ref{ck} below, we give compatibility conditions on these initial polynomials. Namely, we can always  write the initial polynomial $p_0$ as 
$$p_0(\pod \xx)=u_0(\pod \xx)+ v_0(\pod \xx)e_m$$   
for some polynomials $u_0\in\cP^s_k(\bR^{m-1})$ and $v_0\in\cP^{s-1}_k(\bR^{m-1}).$ 
Then we show that the initial polynomial $p_0$ satisfies the compatibility conditions $$\pod\pa^+u_0=0\mand\pod\pa^-v_0=0.$$
In Theorem \ref{ck} below, we prove even that, under the action of $O(m-1),$ 
the CK  extension operator
$CK$ is an invariant isomorphism of the module 
$$\cI^s_k=\Ker^s_k \pod\pa^+\oplus (\Ker^{s-1}_k \pod\pa^-)e_m $$
onto the module $\cH^s_k(\bR^m).$ Here 
$\Ker^s_k \pod\pa^{\pm}=\{u\in\cP^s_k(\bR^{m-1})\ |\ \pod\pa^{\pm}u=0\}.$

To realize the branching law \eqref{branch} for the module $\cH^s_k(\bR^m)$ it is now sufficient to have an irreducible decomposition
$$\cI^s_k=\bigoplus_{\mu_{m-1}}\;\cI(\mu_{m-1})$$
of the module of initial polynomials under the action of $O(m-1).$ Indeed, we then have that
$$\cH^s_k(\bR^m)=\bigoplus_{\mu_{m-1}}\; CK(\cI(\mu_{m-1}))$$
gives a~realization of the branching law \eqref{branch} for the module $\cH^s_k(\bR^m).$ Finally, in Theorem \ref{kerpm} below, we describe irreducible decompositions of the $O(m-1)$-modules $\Ker^s_k \pod\pa^{\pm}$ and thus also $\cI^s_k.$

\paragraph{The Cauchy-Kovalevskaya extension}

Now we characterize restrictions of polynomials of the spaces $\cH^s_k(\bR^m)$ to the hyperplane $x_m=0$ in $\bR^m.$

\begin{thm}\label{ck}
(i) The Cauchy-Kovalevskaya extension operator
$$CK=e^{e_m x_m\pod\pa}$$ is an isomorphism from the module 
$$\cI^s_k=\Ker^s_k \pod\pa^+\oplus (\Ker^{s-1}_k \pod\pa^-)e_m $$
onto the module $\cH^s_k(\bR^m)$ which intertwines the $H$-action of $Pin(m-1).$ 

\medskip\noindent
(ii) Moreover, let $u_0\in\Ker^s_k \pod\pa^+$ and let $v_0\in\Ker^{s-1}_k \pod\pa^-.$ Then we have that
\begin{equation}\label{eqck}
CK(u_0+v_0 e_m)(\xx)=\sum_{j=0}^k\frac{x_m^j}{j!}\;u_j(\pod \xx)+\sum_{j=0}^k\frac{x_m^j}{j!}\;v_j(\pod \xx)e_m
\end{equation}
where
\begin{equation}
\label{equj}	
u_j=\left\{
\begin{array}{ll}
(\pod\pa^+\pod\pa^-)^t u_0,& j=2t,\medskip\\{}
(-1)^{s-1}(\pod\pa^+\pod\pa^-)^t\pod\pa^+v_0,& j=2t+1,
\end{array}
\right.
\end{equation}
\begin{equation}
\label{eqvj}	
v_j=\left\{
\begin{array}{ll}
(\pod\pa^-\pod\pa^+)^t v_0,& j=2t,\medskip\\{}
(-1)^{s-1}(\pod\pa^-\pod\pa^+)^t\pod\pa^-u_0,& j=2t+1.
\end{array}
\right.
\end{equation}
\end{thm}

\begin{proof} First it is well-known that the operator $CK$ is an isomorphism from the space 
of  $\clifford_m$-valued $k$-homogeneous polynomials in $\bR^{m-1}$ onto
the space $$\cM_k(\bR^m)=\{p\in\cP_k^*(\bR^m)\ |\ \pa p=0\}$$ of $k$-homogeneous monogenic polynomials in $\bR^m$ (see \cite[p.\ 152]{DSS}).
Moreover, the operator $CK$ obviously intertwines  the $H$-action of $Pin(m-1)$
since for each polynomial $p\in\cM_k(\bR^m),$ we have that
$$p(\xx)=(e^{e_m x_m\pod\pa}p_0)(\xx)=\sum_{j=0}^k\frac{x_m^j}{j!}\;(e_m\pod\pa)^j p_0(\pod \xx)\text{\ \ with\ \ }p_0(\pod \xx)=p(\pod \xx,0).$$
Now it only remains to show that $CK(\cI^s_k)=\cH^s_k(\bR^m).$

Let $p\in\cH^s_k(\bR^m)$ and let $p_0(\pod \xx)=p(\pod \xx,0).$ We prove that the initial polynomial $p_0$ belongs to the space $\cI^s_k.$ For each $j=1,\ldots,k,$ put $$p_j(\pod \xx)=(e_m\pod\pa)^j p_0(\pod \xx).$$
Of course, we can always write 
$$p_j(\pod \xx)=u_j(\pod \xx)+  v_j(\pod \xx)e_m$$   
for some polynomials $u_j\in\cP^s_{k-j}(\bR^{m-1})$ and $v_j\in\cP^{s-1}_{k-j}(\bR^{m-1}).$ 
Then, for each $j=1,\ldots,k,$ we show that
\begin{equation}\label{equj-1}
u_j=(-1)^{s-1}\pod\pa^+v_{j-1},\ v_j=(-1)^{s-1}\pod\pa^-u_{j-1},\ \pod\pa^+u_{j-1}=0,\ \pod\pa^-v_{j-1}=0.
\end{equation}
Indeed, we have that 
$$
p_j=e_m\pod\pa\; p_{j-1}=e_m\pod\pa^-u_{j-1}+(-1)^{s-1}\pod\pa^+v_{j-1}\mand\pod\pa^+u_{j-1}=0,\ \pod\pa^-v_{j-1}=0
$$
and this since $\pod\pa^+u_{j-1}\in\cP^{s+1}_{k-j}(\bR^{m-1})$ and $\pod\pa^-v_{j-1}\in\cP^{s-2}_{k-j}(\bR^{m-1}).$ 
Hence, by \eqref{equj-1}, we easily get that $p_0\in\cI^s_k$ and the formul\ae\ \eqref{equj} and \eqref{eqvj}.

On the other hand, for a~given $p_0\in\cI^s_k,$ we show that the polynomial $CK(p_0)$ belongs to the space $\cH^s_k(\bR^m).$ But, in this case, it is easy to get the formul\ae\ \eqref{eqck}, \eqref{equj} and \eqref{eqvj}, 
which finishes the proof.
\end{proof}

\paragraph{The Fischer decompositions of the spaces $\Ker^s_k \pa^{\pm}$}

Let us first define the Euler operator $E$ and the fermionic Euler operators $\pa^+\rfloor$ and $\pa^-\rceil$ by
\begin{equation}\label{eulers}
E=\sum_{j=1}^mx_j\pa_{x_j},\ \ \ \pa^+\rfloor=-\sum_{j=1}^m (e_j\wed) (e_j\bul)\text{\ \ \ and\ \ \ }
\pa^-\rceil=-\sum_{j=1}^m (e_j\bul) (e_j\wed).
\end{equation}
If $P\in\cP^s_k,$ then it is easy to see that $$EP=kP,\ \ \ \pa^+\rfloor P=sP\text{\ \ and\ \ }\pa^-\rceil P=(m-s)P$$
(see \cite{BDS} for details).
Putting $A=E+\pa^+\rfloor$ and $B=E+\pa^-\rceil,$ we have that for each $P\in\cP^s_k,$ $AP=(s+k)P$ and $BP=(m-s+k)P.$ 
Furthermore, in a classical way the Laplace operator  $\Delta$ in $\bR^m$ is defined by
$$\Delta=\sum_{j=1}^m {\pa_{x_j}^2}.$$
Now we are ready to describe 
the Fischer decompositions of the spaces $\Ker^s_k \pa^{\pm}.$ 

\begin{thm}\label{kerpm} Let $1\leq s\leq m-1.$ Then the following statements hold:

\medskip\noindent
(i) Under the $H$-action, the space $\Ker^s_k\pa^+$ has the multiplicity free irreducible decomposition
\begin{equation}\label{eqkerp}
\Ker^s_k\pa^+=\cH^s_k\oplus\bigoplus_{j=0}^{[(k-1)/2]}\xx^{2j}(\xx\wed)\cH^{s-1}_{k-2j-1}
\oplus\bigoplus_{j=0}^{[(k-2)/2]} \dot\yy_{2j+2}\cH^s_{k-2j-2}
\end{equation}
where $\dot\yy_{2j+2}=\xx^{2j+1}((\xx\bul)(A+2j+2)+(\xx\wed) A)$ with $A=E+\pa^+\rfloor.$ 

\medskip\noindent
(ii) Under the $H$-action, the space $\Ker^s_k\pa^-$ has the multiplicity free irreducible decomposition
\begin{equation}\label{eqkerm}
\Ker^s_k\pa^-=\cH^s_k\oplus\bigoplus_{j=0}^{[(k-1)/2]}\xx^{2j}(\xx\bul)\cH^{s+1}_{k-2j-1}
\oplus\bigoplus_{j=0}^{[(k-2)/2]}\hat\yy_{2j+2}\cH^s_{k-2j-2}
\end{equation}
where $\hat\yy_{2j+2}=\xx^{2j+1}((\xx\wed)(B+2j+2)+(\xx\bul) B)$ with $B=E+\pa^-\rceil.$

\end{thm}

\begin{rem} 
It is easy to see that  
$\Ker^0_k\pa^-=\cP^0_k,$ $\Ker^0_k\pa^+=\cH^0_k,$ $\Ker^m_k\pa^+=\cP^m_k$ and $\Ker^m_k\pa^-=\cH^m_k.$
\end{rem}

Before proving Theorem \ref{kerpm} we need some lemmas. 

\begin{lem}\label{lrels}
If   for linear operators $T$ and $S$ acting on the space $\cP^*$ we put $\{T,S\}=TS+ST,$ then we have that
\begin{equation*}
\begin{array}{lll}
\{(\xx\wed),(\xx\wed)\}=0, &\{(\xx\bul),(\xx\bul)\}=0, &\{(\xx\wed),(\xx\bul)\}=\xx^2,\medskip\\{}
\{\pa^+,\pa^+\}=0, &\{\pa^-,\pa^-\}=0, &\{\pa^+,\pa^-\}=-\Delta,\medskip\\{}
\{(\xx\bul),\pa^+\}=-A, &\{(\xx\wed),\pa^-\}=-B, &\{(\xx\bul),\pa^-\}=0=\{(\xx\wed),\pa^+\}.
\end{array}
\end{equation*}
\end{lem}

\begin{proof}See e.g. \cite{BDS}.\end{proof}

Using Lemma \ref{lrels}, it is easy to prove the next relations.

\begin{lem}\label{lrelsj}
If  for linear operators $T$ and $S$ acting on the space $\cP^*$ we put $[T,S]=TS-ST,$  then we have that
\begin{equation*}
\begin{array}{ll}
[\pa^+,\xx^{2j+1}(\xx\bul)]=\xx^{2j}(\xx\wed) A,& 
[\pa^+,\xx^{2j+1}(\xx\wed)]=-\xx^{2j}(\xx\wed) (A+2j+2), \medskip\\{}  
[\pa^-,\xx^{2j+1}(\xx\wed)]=\xx^{2j}(\xx\bul) B,& 
[\pa^-,\xx^{2j+1}(\xx\bul)]=-\xx^{2j}(\xx\bul) (B+2j+2), \medskip\\{}  
[\pa^+,\xx^{2j}]=-2j\;\xx^{2(j-1)}(\xx\wed),&
[\pa^-,\xx^{2j}]=-2j\;\xx^{2(j-1)}(\xx\bul), \medskip\\{}
\{\pa^+, \xx^{2j+2}(\xx\bul)\}=\dot\yy_{2j+2}, &
\{\pa^-, \xx^{2j+2}(\xx\wed)\}=\hat\yy_{2j+2}.
\end{array}
\end{equation*}
In addition, on the space $\cH^s_{k-2j-2},$ we have that
$$\pa^-\dot\yy_{2j+2}=-(2j+2)(m+2k-2j-2)\;\xx^{2j}(\xx\bul),\ \ \pa^+\dot\yy_{2j+2}=0;$$
$$\pa^+\hat\yy_{2j+2}=-(2j+2)(m+2k-2j-2)\;\xx^{2j}(\xx\wed),\ \ \pa^-\hat\yy_{2j+2}=0.$$
\end{lem}

In the~proof of the Fischer decompositions of the spaces $\Ker^s_k\pa^{\pm},$ we shall also use the next decompositions.

\begin{prop}\label{lfp}
We have that
$$\cP^s_k=\Ker^s_k\pa^+\oplus\; (\xx\bul)\Ker^{s+1}_{k-1}\pa^+.$$ 
Moreover, the projection $P^+$ of the space $\cP^s_k$ onto the space $\Ker^s_k\pa^+$ is given by
$$P^+=-(s+k)^{-1}\pa^+(\xx\bul).$$
\end{prop}

\begin{proof}
Obviously, using the relation $\{(\xx\bul),\pa^+\}=-A$ of Lemma \ref{lrels}, we have that 
$$\Ker^s_k\pa^+\cap\;(\xx\bul)\Ker^{s+1}_{k-1}\pa^+=\{0\}.$$ Furthermore, if
for a~given polynomial $p\in\cP^s_k,$ we put 
$$p^+=-\pa^+(\xx\bul) A^{-1} p\mand p^-=-(\xx\bul)\pa^+A^{-1} p,$$
then it is easily seen that $p=p^++p^-$ with $p^+\in\Ker^s_k\pa^+$ and $p^-\in(\xx\bul)\Ker^{s+1}_{k-1}\pa^+.$
This completes the proof.
\end{proof}

Of course, we can prove an analogous proposition for the operator $\pa^-.$

\begin{prop}\label{lfm}
We have that
$$\cP^s_k=\Ker^s_k\pa^-\oplus (\xx\wed)\Ker^{s-1}_{k-1}\pa^-.$$ 
Moreover, the projection $P^-$ of the space $\cP^s_k$ onto the space $\Ker^s_k\pa^-$ is given by
$$P^-=-(m-s+k)^{-1}\pa^-(\xx\wed).$$
\end{prop}

\begin{proof}[Proof of Theorem \ref{kerpm}]
Using the Fischer decomposition for the $H$-action (see Theorem \ref{tFischerH}), we get the next irreducible (not multiplicity free) decomposition of the space $\cP^s_k:$
$$\cP^s_k=\cH^s_k\oplus\bigoplus_{j=0}^{[(k-1)/2]}\xx^{2j}(\xx\wed)\cH^{s-1}_{k-2j-1}\oplus
\bigoplus_{j=0}^{[(k-1)/2]}\xx^{2j}(\xx\bul)\cH^{s+1}_{k-2j-1}\oplus$$
$$\oplus\bigoplus_{j=0}^{[(k-2)/2]}\xx^{2j+1}(\xx\bul)\cH^s_{k-2j-2}
\oplus\bigoplus_{j=0}^{[(k-2)/2]}\xx^{2j+1}(\xx\wed)\cH^s_{k-2j-2}.$$
Applying the projections $P^{\pm}$ of Propositions \ref{lfp} and \ref{lfm} to this decomposition, we easily get the required decompositions of the spaces $\Ker^s_k\pa^{\pm}.$ 
Indeed, by Lemmas \ref{lrels} and \ref{lrelsj}, we have that $P^{\pm}(\cP^s_k)=\Ker^s_k\pa^{\pm},$ $P^{\pm}(\cH^s_k)=\cH^s_k,$ 
\begin{equation*}
\begin{array}{lll}
P^+(\xx^{2j}(\xx\wed)\cH^{s-1}_{k-2j-1}) &=& \xx^{2j}(\xx\wed)\cH^{s-1}_{k-2j-1},\medskip\\{}
P^+(\xx^{2j+1}(\xx\bul)\cH^s_{k-2j-2}) &=& \dot\yy_{2j+2}\cH^s_{k-2j-2},\medskip\\{} 
P^-(\xx^{2j}(\xx\bul)\cH^{s+1}_{k-2j-1}) &=& \xx^{2j}(\xx\bul)\cH^{s+1}_{k-2j-1},\medskip\\{}
P^-(\xx^{2j+1}(\xx\wed)\cH^s_{k-2j-2}) &=& \hat\yy_{2j+2}\cH^s_{k-2j-2}.
\end{array}
\end{equation*}
Moreover, the projections $P^{\pm}$ vanish on the remaining pieces, which finishes the proof.
\end{proof}

\section{Explicit description of GT bases}

In this section, we use the CK method explained in the previous section to construct quite explicitly GT bases for the spaces $\cH^s_k$ of solutions of Hodge-de Rham systems. 

\paragraph{Induction step}

First we explain how to construct GT bases for solutions  of the Hodge-de Rham systems in $\bR^m$ when we already know these bases in $\bR^{m-1}.$ To do this we need some lemmas. But first recall that the Gegenbauer polynomial $C^{\nu}_j$ is defined as
\begin{equation}\label{gegenbauer}
C^{\nu}_j(z)=\sum_{i=0}^{[j/2]}\frac{(-1)^i(\nu)_{j-i}}{i!(j-2i)!}(2z)^{j-2i}\text{\ \ with\ \ }
(\nu)_{j}=\nu (\nu+1)\cdots (\nu+j-1).
\end{equation}

\begin{lem}\label{lckxj} 
Let $j\in\bN_0$ and let $P_k\in\cM_k(\bR^{m-1}).$ Then we have that
$$CK(\pod\xx^jP_k(\pod\xx))=X^j_k(\pod\xx,x_m)P_k(\pod\xx)$$
where $X^{0}_k=1$ and, for $j\in\bN,$ the polynomial $X^{j}_k$ is given by
$$
X^j_k(\pod\xx,x_m)=
\mu^j_kr^j\left(C_j^{m/2+k-1}(\frac{x_m}{r})+
\frac{m+2k-2}{m+2k+j-2}C_{j-1}^{m/2+k}(\frac{x_m}{r})\frac{\nad e_m\pod\xx}{r}\right)\;\nad e_m^j
$$ 
with $r=(x_1^2+x^2_2+\cdots+x_m^2)^{1/2},$
$\mu^{2l}_k=(C_{2l}^{m/2+k-1}(0))^{-1}$ and 
$$\mu^{2l+1}_k=\frac{m+2k+2l-1}{m+2k-2}(C_{2l}^{m/2+k}(0))^{-1}.$$
\end{lem}

\begin{proof}
In \cite[p. 312, Theorem 2.2.1]{DSS}, the corresponding polynomial we denote here by $\tilde X^j_k$ is computed for the Cauchy-Riemann operator. Fortunately,
there is an obvious relation between these two polynomials. Indeed, we have that
\begin{equation*}
X^j_k(\pod\xx,x_m)=\left\{
\begin{array}{ll}
\tilde X^j_k(\nad e_m\pod\xx,x_m),& j\text{\ even},\medskip\\{}
\tilde X^j_k(\nad e_m\pod\xx,x_m)\nad e_m,& j\text{\ odd}.
\end{array}
\right.
\end{equation*}
To complete the proof it is sufficient to use the explicit formula for the polynomial $\tilde X^j_k.$
\end{proof}

\begin{lem}\label{lckxjo} 
Let either $P\in\cH^s_k(\bR^{m-1})$ or $P= Q e_m$ for some $Q\in\cH^s_k(\bR^{m-1}).$ 

\medskip\noindent
(i) We have that
$CK(\pod\xx^{j-1}(\pod\xx\wed) P(\pod\xx))=\hat X_j P(\pod\xx)$
where $$\hat X_j=\hat X_j^{s,k}=X^{j-1}_{k+1}(\pod\xx\wed)+(1-c)\;(X^j_k-X^{j-1}_{k+1}\pod\xx)$$ with $c=(s+k)(m-1+2k)^{-1}.$

\medskip\noindent
(ii) We have that
$CK(\pod\xx^{j-1}(\pod\xx\bul) P(\pod\xx))=\dot X_j P(\pod\xx)$
where $$\dot X_j=\dot X_j^{s,k}=X^{j-1}_{k+1}(\pod\xx\bul)+c\;(X^j_k-X^{j-1}_{k+1}\pod\xx)$$ 
with the constant $c$ being the same as in (i).

\medskip\noindent
(iii) We have that
$CK(\dot\pyy_{2j+2} P(\pod\xx))=\dot Y_{2j+2} P(\pod\xx)$
where 
$$\dot Y_{2j+2}=\dot Y_{2j+2}^{s,k}=(s+k)X^{2j+2}_k+(2j+2)\dot X_{2j+2}.$$
Here $\dot\pyy_{2j+2}=\pod\xx^{2j+1}((s+k+2j+2)(\pod\xx\bul)+(s+k)(\pod\xx\wed)).$

\medskip\noindent
(iv) We have that
$CK(\hat\pyy_{2j+2} P(\pod\xx))=\hat Y_{2j+2} P(\pod\xx)$
where 
$$\hat Y_{2j+2}=\hat Y_{2j+2}^{s,k}=(m-1-s+k)X^{2j+2}_k+(2j+2)\hat X_{2j+2}.$$
Here $\hat\pyy_{2j+2}=\pod\xx^{2j+1}((m+1-s+k+2j)(\pod\xx\wed)+(m-1-s+k)(\pod\xx\bul)).$

\medskip\noindent
In addition, we have that $\hat X_j^{s,k}+\dot X_j^{s,k}=X^j_k$ and that
$$\hat Y_{2j+2}^{s,k}+\dot Y_{2j+2}^{s,k}=(m+1+2k+2j)\,X^{2j+2}_k.$$
\end{lem}

\begin{proof}
We prove the formula in (i).
Obviously, the polynomial $H(\pod\xx)=\pod\xx\wedge P(\pod\xx)$ is harmonic
and, as is well-known, we then have that $H(\pod\xx)=M_0(\pod\xx)+\pod\xx M_1(\pod\xx)$ 
for some monogenic polynomials $M_0$ and $M_1.$ It may be easily checked that 
$$M_1(\pod\xx)=(1-c) P(\pod\xx)\mand 
M_0(\pod\xx)=((\pod\xx\wed)-(1-c)\;\pod\xx) P(\pod\xx).$$
Moreover, by Lemma \ref{lckxj}, we get that 
$$CK(\pod\xx^{j-1}H(\pod\xx))=CK(\pod\xx^{j-1}M_0(\pod\xx))+CK(\pod\xx^{j}M_1(\pod\xx))
=X^{j-1}_{k+1}M_0(\pod\xx)+X^{j}_{k}M_1(\pod\xx),$$
which easily completes  the proof.

Of course, we can show the formula in (ii) in quite an analogous way.
The remaining relations in (iii) and (iv) are then obvious.
\end{proof}

Now we are ready to prove the following theorem.

\begin{thm}\label{induction}
Let $1\leq s\leq m-1$ and let $k\in\bN_0.$ Furthermore, let for each $t= s-1, s$ and $l=0,\ldots, k,$  $\cB^{t,m-1}_{l}$
stand for a~GT basis of the space $\cH^{t}_{l}(\bR^{m-1}).$ Then the space $\cH^s_k(\bR^m)$ has  a~GT basis $\cB^{s,m}_k=\cB^+\cup\cB^-$ where
$$
\cB^+=\cB^{s,m-1}_k\cup\bigcup_{j=0}^{[(k-1)/2]}\hat X_{2j+1}\cB^{s-1,m-1}_{k-2j-1}
\cup\bigcup_{j=0}^{[(k-2)/2]} \dot Y_{2j+2}\cB^{s,m-1}_{k-2j-2}
\text{\ \ \ and\ \ \ }
$$
$$
\cB^-=\cB^{s-1,m-1}_k e_m\cup\bigcup_{j=0}^{[(k-1)/2]}\dot X_{2j+1}\cB^{s,m-1}_{k-2j-1}e_m
\cup\bigcup_{j=0}^{[(k-2)/2]} \hat Y_{2j+2}\cB^{s-1,m-1}_{k-2j-2}e_m.
$$
Here we denote, for example, $\hat X_{2j+1}\cB^{s-1,m-1}_{k-2j-1}=\{\hat X_{2j+1}P |\ P\in \cB^{s-1,m-1}_{k-2j-1}\}$ and 
$\cB^{s-1,m-1}_k e_m=\{P e_m|\ P\in \cB^{s-1,m-1}_k\}.$
\end{thm}

\begin{proof} By Theorem \ref{ck}, 
we know that the CK extension operator is an invariant isomorphism from the space
$$\cI^s_k=\Ker^s_k \pod\pa^+\oplus (\Ker^{s-1}_k \pod\pa^-)e_m $$
of initial polynomials onto the space $\cH^s_k(\bR^m).$ 
Moreover, Theorem \ref{kerpm} tells us that the space $\cI^s_k$ has a~basis
$b^s_k=b^+\cup b^-$ where  
$$
b^+=\cB^{s,m-1}_k\cup\bigcup_{j=0}^{[(k-1)/2]}\pod\xx^{2j}(\pod\xx\wed)\,\cB^{s-1,m-1}_{k-2j-1}
\cup\bigcup_{j=0}^{[(k-2)/2]} \dot\pyy_{2j+2}\cB^{s,m-1}_{k-2j-2}
\text{\ \ \ and\ \ \ }
$$
$$
b^-=\cB^{s-1,m-1}_k e_m\cup\bigcup_{j=0}^{[(k-1)/2]}\pod\xx^{2j}(\pod\xx\bul)\, \cB^{s,m-1}_{k-2j-1}e_m
\cup\bigcup_{j=0}^{[(k-2)/2]} \hat\pyy_{2j+2} \cB^{s-1,m-1}_{k-2j-2} e_m.
$$
As we explained before, we get the GT basis $\cB^{s,m}_k$ for the space $\cH^s_k(\bR^m)$ by applying the CK extension operator to the elements of the basis $b^s_k,$ i.e., $CK(b^s_k)=\cB^{s,m}_k.$ 
To finish the proof it is now sufficient to use Lemma \ref{lckxjo}.
\end{proof}

\paragraph{Examples} For $\clifford_m=\bC_m$ (resp.\ $\clifford_m=\bR_{0,m}$), we describe below GT bases $\cB^{s,m}_k$
(resp.\ $\tilde\cB^{s,m}_k$) of the spaces $\cH^s_k(\bR^m)$ of solutions of the Hodge-de Rham system in some special cases. If we construct GT bases for a~particular space $\cH^s_k(\bR^m)$ in various ways we do not claim that all of these bases are identical. On the other hand, we know,  by the definition,  that the corresponding elements of these bases must be the same up to non-zero multiples. 

\begin{ex}
We can put
$\cB^{0,m}_0=\{1\},$ 
$\cB^{m,m}_0=\{e_M\}$ with $e_M=e_1e_2\cdots e_m$ and, for $s\in\{0,m\}$ and $k>0,$ we have that $\cB^{s,m}_k=\emptyset.$
\end{ex}

\begin{ex}
{\it The Riesz system} (i.e., the case when $s=1$): Assume first that $\clifford_m=\bC_m.$
It is well-known (see \cite[p.\ 460]{vil} for more details) that a~canonical basis (i.e., a~GT basis) $B^m_{k+1}$ of the space $\Harm_{k+1}(\bR^m)$ of 
complex valued 
spherical harmonics of degree $k+1$ in $\bR^m$  is formed (up to normalization) by the polynomials 
\begin{equation}\label{harmonics}	
\Xi^{k+1}_{\mu,\pm}=(x_1\pm i x_2)^{k_{m-2}}
\prod_{j=0}^{m-3}r_{m-j}^{k_j-k_{j+1}}C^{(m-j-2)/2+k_{j+1}}_{k_j-k_{j+1}}\left(\frac{x_{m-j}}{r_{m-j}}\right)
\end{equation}
where $r^2_{m-j}=x_1^2+\cdots+x^2_{m-j}$ and  
$\mu$ is an arbitrary sequence of integers $(k_1,\ldots, k_{m-2})$ such that $k+1=k_0\geq k_1\geq\ldots\geq k_{m-2}\geq 0.$ 
Then we can put $$\cB^{1,m}_k=\{\pa\; \Xi\ |\ \Xi\in B^m_{k+1} \}.$$ Indeed, it is easy to see that the Dirac operator $\pa$ (or also $\pa^+$) is an invariant isomorphism from the $O(m)$-module $\Harm_{k+1}(\bR^m)$ onto the module $\cH^1_k(\bR^m).$  
Actually, these bases are well-known  (see \cite{cno}). 
Furthermore, we can take
$$
\cB^{m-1,m}_k=\{\pa\; \Xi\;e_M|\ \Xi\in B^m_{k+1} \}.
$$
Moreover $\cB^{s,m}_ke_M$ is a~GT basis of $\cH^{m-s}_k(\bR^m)$
whenever $\cB^{s,m}_k$ is a~GT basis of $\cH^s_k(\bR^m).$

In the case when $\clifford_m=\bR_{0,m},$ we construct the GT bases $\tilde\cB^{1,m}_k$ and $\tilde\cB^{m-1,m}_k$ in quite an analogous way using, of course, in this case the canonical basis
$$\tilde B^m_{k+1}=\{\Re\,\Xi^{k+1}_{\mu,+}, \Im\,\Xi^{k+1}_{\mu,+}  \}$$
of the space of real valued spherical harmonics of degree $k+1.$
Here, $\mu$ is as in (\ref{harmonics}) and for a~complex number $z,$ $\Re z$ is its real part and $\Im z$ is its imaginary part. 
\end{ex}

\begin{ex} \label{exGTB34}
Let $\clifford_m=\bC_m.$ We explain now how to construct GT bases $\cB^{s,m}_k$ using the CK method. 
In the case when the dimension $m=2,$ we know, by Example 1 and Example 2, that
\begin{equation}\label{GT2}
\cB^{0,2}_0=\{1\},\ \ \cB^{2,2}_0=\{e_{12}\}\text{\ \ and\ \ }\cB^{1,2}_k=\{z_{\pm}^k w_{\pm}\}
\end{equation}
where $z_{\pm}=x_1\pm ix_2$ and $w_{\pm}=e_1\pm ie_2.$
Obviously, we can obtain GT bases in higher dimensions inductively with the help of Theorem \ref{induction}.
For example, in dimension 3 we have that
\begin{eqnarray*}
\cB^{1,3}_0 &=& \{w_{\pm},\ e_3\},\\ 
\cB^{1,3}_1 &=& \{z_{\pm}w_{\pm},\ (z_-w_+ +z_+w_-)/2 -2x_3e_3,\  -z_{\pm}e_3-x_3w_{\pm}\},\\
\cB^{1,3}_2 &=& \{z_{\pm}^2w_{\pm},\ -2x_3z_{\pm}w_{\pm}-z_{\pm}^2e_3,\\ 
&&\ \ (8x_3^2-4z_+z_-)e_3-4x_3(z_-w_+ +z_+w_-),\\ 
&&\ \ 8x_3z_{\pm}e_3-z_{\pm}^2w_{\mp}+(4x_3^2-2z_+z_-)w_{\pm}
\}.
\end{eqnarray*}
Moreover, we know that $\cB^{2,3}_k=\cB^{1,3}_k e_{123}.$  

By Theorem \ref{induction}, we can, for example, compute the following GT bases of bivector valued monogenic polynomials in dimension 4: 

\medskip\noindent
$\cB^{2,4}_0 = \{e_{12},\, e_{34},\, w_{\pm}e_3,\, w_{\pm}e_4\},$

\medskip\noindent
$\cB^{2,4}_1 = \cB^{2,3}_1\,\cup\,\cB^{1,3}_1e_4\,\cup\,$\\ \smallskip
$\mbox{\ \ \ \ \ \ \ \ \ }
\{(z_+w_- +z_-w_+)e_3/2+ 2x_4e_{34},\,
-x_3w_{\pm}e_3\pm i z_{\pm}e_{12}+2x_4w_{\pm}e_4,$\\ \smallskip
$\mbox{\ \ \ \ \ \ \ \ \ }
i(z_-w_+ -z_+w_-)e_4/2+2x_4e_{12},\, - z_{\pm}e_{34}+ x_3w_{\pm}e_4+ 2 x_4w_{\pm} e_3\},$

\medskip\noindent
$\cB^{2,4}_2 = \cB^{2,3}_2\, \cup\, \cB^{1,3}_2e_4\, \cup\,$\\ \smallskip
$\mbox{\ \ \ \ \ \ \ \ \ } \{-x_3z_{\pm}w_{\pm}e_3\pm z_{\pm}^2 i e_{12}+3x_4z_{\pm}w_{\pm}e_4,$\\ \smallskip
$\mbox{\ \ \ \ \ \ \ \ \ }
(3/2)(z_-w_+ +z_+w_-)(x_4e_4-x_3e_3)-6x_3x_4e_{34},$\\ \smallskip
$\mbox{\ \ \ \ \ \ \ \ \ }
x_3^2w_{\pm}e_3-(z_{\pm}/2)(z_-w_+ +z_+w_-)e_3\mp x_3z_{\pm}i e_{12}-3x_4z_{\pm}e_{34}-3x_3x_4w_{\pm}e_4,$\\ \smallskip
$\mbox{\ \ \ \ \ \ \ \ \ }
i(z_-w_+ -z_+w_-)(5x_4e_4-x_3e_3)-(4z_+z_- +2x_3^2-10x_4^2)e_{12},$\\ \smallskip
$\mbox{\ \ \ \ \ \ \ \ \ }
(10x_4^2-4x_3^2-3z_+z_-)w_{\pm}e_3 + 10x_3x_4w_{\pm}e_4 -z_{\pm}^2w_{\mp}e_3\pm 2x_3z_{\pm}ie_{12}-10x_4z_{\pm}e_{34},$\\ \smallskip
$\mbox{\ \ \ \ \ \ \ \ \ }
x_3^2w_{\pm}e_4\mp (z_{\pm}/2)(z_-w_+ -z_+w_-)e_4 - x_3z_{\pm}e_{34}\pm 3x_4z_{\pm}ie_{12}+3x_3x_4w_{\pm}e_3,$\\ \smallskip
$\mbox{\ \ \ \ \ \ \ \ \ } 
x_3z_{\pm}w_{\pm}e_4- z_{\pm}^2 e_{34}+3x_4z_{\pm}w_{\pm}e_3,$\\ \smallskip
$\mbox{\ \ \ \ \ \ \ \ \ }
(3/2)i(z_-w_+ -z_+w_-)(x_3e_4+x_4e_3)+6x_3x_4e_{12},$\\ \smallskip
$\mbox{\ \ \ \ \ \ \ \ \ }
(10x_4^2-4x_3^2-3z_+z_-)w_{\pm}e_4 - 10x_3x_4w_{\pm}e_3 +z_{\pm}^2w_{\mp}e_4 + 2x_3z_{\pm}e_{34} \pm  10x_4z_{\pm}i e_{12},$\\ \smallskip
$\mbox{\ \ \ \ \ \ \ \ \ }
(z_-w_+ +z_+w_-)(5x_4e_3+x_3e_4)-(4z_+z_- +2x_3^2-10x_4^2)e_{34}\}.$

To summarize we have an algorithm how to obtain any particular GT basis $\cB^{s,m}_k$ by induction on the dimension $m.$
Actually, all explicit examples in this paper were computed using the mathematical software Maple and the Maple package Clifford (see \cite{AF}).  

Now it remains to deal with the case when $\clifford_m=\bR_{0,m}.$ Let us notice that all the polynomials $\hat X_{2j+1},$ $\dot X_{2j+1},$ $\hat Y_{2j+2}$ and $\dot Y_{2j+2}$ are $\bR_{0,m}$-valued.  When we thus start  with the GT bases \eqref{GT2} in dimension 2, the explained construction gives us GT bases $\cB^{s,m}_k$ with the following property: Each basis $\cB^{s,m}_k$ consists partly of $\bR_{0,m}$-valued basis elements $P_{\alpha}$ and partly of pairs of complex conjugate basis elements  $P^{\pm}_{\beta}.$ Obviously, we can make a~'real' GT basis $\tilde\cB^{s,m}_k$ from the 'complex' basis $\cB^{s,m}_k$ just by replacing each pair of complex conjugate basis elements  $P^{\pm}_{\beta}$ of the basis $\cB^{s,m}_k$
with the pair $\Re\, P^+_{\beta}$ and $\Im\, P^+_{\beta}.$      

\end{ex}

\paragraph{Generalized Moisil-Th\'eodoresco  systems} Now we construct orthogonal bases for spaces of homogeneous solutions of GMT systems. 
Let $S$ be a~subset of $\{0,1,\ldots,m\}.$ 
Recall that $\cM^S_k(\bR^m)$ stands for 
the space of $k$-homogeneous $\clifford_m^S$-valued spherical monogenics in $\bR^m,$
where $$\clifford_m^S=\bigoplus_{s\in S}\clifford_m^s.$$
The following decomposition of this space is known (see \cite{DLS2, lav_isaac09}).

\begin{thm}\label{tGMT} Let $S\subset\{0,1,\ldots,m\}$ and $S'=\{s:s\pm 1\in S\}.$
Under the $H$-action, the space $\cM^S_k(\bR^m)$ decomposes into inequivalent irreducible pieces as
$$\cM^S_k(\bR^m)=\left(\bigoplus_{s\in S} \cH^s_k\right)\oplus\left(\bigoplus_{s\in S'}((k-1+m-s)(\xx\bul)-(k-1+s)(\xx\wed))\,\cH^s_{k-1}\right).$$
\end{thm}

A~direct consequence of Theorem \ref{tGMT} is

\begin{cor}\label{OG_GMT}
The space $\cM^S_k(\bR^m)$ has a~basis
$$\cB_k^{S,m}=\left(\bigcup_{s\in S} \cB^{s,m}_k\right)\cup\left(\bigcup_{s\in S'}((k-1+m-s)(\xx\bul)-(k-1+s)(\xx\wed))\,\cB^{s,m}_{k-1}\right).$$
Here $\cB^{s,m}_k$ is a~GT basis of the space $\cH^s_k(\bR^m).$
In particular, the space $\cM_k(\bR^m)$ of $k$-homogeneous $\clifford_m$-valued spherical monogenics in $\bR^m$ has an orthogonal basis $\cB^{M,m}_k$ with $M=\{0,1,\ldots,m\}.$ 

Moreover, the basis $\cB_k^{S,m}$ is orthogonal with respect to any invariant inner product, including the $L_2$-inner product \eqref{L2product} and the Fischer inner product \eqref{Fischerproduct}.
\end{cor}

\begin{ex} Assume that $\clifford_m=\bC_m.$
According to Corollary \ref{OG_GMT}, an orthogonal basis $\cB_k^{S,m}$ for the corresponding GMT system includes partly GT bases $\cB^{s,m}_k$ of spaces $\cH^s_k(\bR^m)$ for $s\in S$ and partly subsets of the form
$$\cV^{s,m}_k=((k-1+m-s)(\xx\bul)-(k-1+s)(\xx\wed))\,\cB^{s,m}_{k-1},\ \ s\in S'.$$
Of course, here $\cV^{s,m}_0=\emptyset.$ Moreover, we have that (up to a~normalization)  $\cV^{m-s,m}_k=\cV^{s,m}_k e_M.$ Hence, using Example \ref{exGTB34}, to describe explicitly the orthogonal bases for all GMT systems in dimension 3 with $k=0,1,2$ it is sufficient to compute just the following sets:

\medskip\noindent
$\cV^{1,3}_1=\{-2z_{\pm}+x_3w_{\pm}e_3\mp z_{\pm}i e_{12},\ -2x_3-(z_-w_+ +z_+w_-)e_3/2\},$

\medskip\noindent
$\cV^{1,3}_2=\{-3z_{\pm}^2+2x_3z_{\pm}w_{\pm}e_3\mp 2z_{\pm}^2i e_{12},\
6x_3^2-3z_+z_- +3(z_-w_+ +z_+w_-)x_3e_3,$\\ \smallskip
$\mbox{\ \ \ \ \ \ \ \ \ \ \ }
6x_3z_{\pm}\pm 2x_3z_{\pm}ie_{12}+(z_+z_- -2x_3^2)w_{\pm}e_3+z_{\pm}^2w_{\mp}e_3\}.$

\medskip\noindent
Actually, we have again an algorithm how to obtain any particular orthogonal basis $\cB^{S,m}_k.$
Furthermore, the case when $\clifford_m=\bR_{0,m}$ can be dealt with as in Example \ref{exGTB34}. 
\end{ex}

\subsection*{Acknowledgments}

R. L\'avi\v cka and V. Sou\v cek acknowledge the financial support from the grant GA 201/08/0397.
This work is also a part of the research plan MSM 0021620839, which is financed by the Ministry of Education of the Czech Republic.



\bigskip\bigskip\bigskip

\noindent
Richard Delanghe,\\ Clifford Research Group, Department of Mathematical Analysis,\\ Ghent University, Galglaan 2, B-9000 Gent, Belgium\\
email: \texttt{richard.delanghe@ugent.be}

\bigskip

\noindent
Roman L\'avi\v cka and Vladim\'ir Sou\v cek,\\ Mathematical Institute, Charles University,\\ Sokolovsk\'a 83, 186 75 Praha 8, Czech Republic\\
email: \texttt{lavicka@karlin.mff.cuni.cz} and \texttt{soucek@karlin.mff.cuni.cz}


\begin{thebibliography}{1}
%
\bibitem{AF} R. Ab{\l}amowicz and  B. Fauser, \textit{CLIFFORD/Bigebra, A Maple Package for Clifford (Co)Algebra Computations}, 2009 (available at http://www.math.tntech.edu/rafal).
%
\bibitem{BGLS} S. Bock, K. G\"urlebeck, R. L\' avi\v cka and V. Sou\v cek, The Gelfand-Tsetlin bases for spherical monogenics in dimension 3, preprint.

 \bibitem{Bock2010c} S. Bock, Orthogonal Appell bases in dimension 2,3 and 4. In Numerical Analysis and Applied Mathematics
(T.E. Simos, G. Psihoyios, and Ch. Tsitouras, eds.), AIP Conference Proceedings, vol. 1281. American
Institute of Physics: Melville, NY, 2010; 1447--1450.
%
\bibitem{Bock2010a} S. Bock, On a three dimensional analogue to the holomorphic $z$-powers: Power series and recurrence
formulae, submitted, 2010.
%
\bibitem{Bock2009} S. Bock, \"{U}ber funktionentheoretische {M}ethoden in der r\"{a}umlichen
{E}lastizit\"{a}tstheorie, PhD thesis, Bauhaus-University, Weimar, (url:
http://e-pub.uni-weimar.de/frontdoor.php?source\_opus=1503, date: 07.04.2010), 2009.
%
\bibitem{BG} S. Bock and K. G\"urlebeck, On a~generalized Appell
system and monogenic power series, Mathematical Methods in the Applied Sciences 33 (2010),
394--411.
%
\bibitem{BDS} F. Brackx, R. Delanghe and F. Sommen, \textit{Differential forms and/or multi-vector
functions}, CUBO \textbf{7} (2005), 139-170.
%
\bibitem{BDS2} F. Brackx, R. Delanghe and F. Sommen, \textit{Clifford analysis}, Pitman, London, 1982.

\bibitem{BSES} F. Brackx, H. De Schepper, D. Eelbode and V. Sou\v cek, \textit{The Howe dual pair in hermitian Clifford analysis}, 
Rev. Mat. Iberoamericana \textbf{26} (2010)(2), 449--479.
%
\bibitem{ckH} F.\ Brackx, H.\ De Schepper, R.\ L\'{a}vi\v{c}ka, V.\ Sou\v{c}ek,
\textit{The Cauchy-Kovalevskaya	Extension Theorem in Hermitean Clifford Analysis}, preprint.
%
\bibitem{kerH} F.\ Brackx, H.\ De Schepper, R.\ L\'{a}vi\v{c}ka, V.\ Sou\v{c}ek, 
\textit{Fischer decompositions of kernels of Hermitean Dirac operators}, In: T.E.\ Simos, G.\ Psihoyios, Ch.\ Tsitouras, {\em Numerical Analysis and Applied Mathematics}, AIP Conference Proceedings, Rhodes, Greece (2010).
%
\bibitem{GTinH} F.\ Brackx, H.\ De Schepper, R.\ L\'{a}vi\v{c}ka, V.\ Sou\v{c}ek, 
\textit{Gel'fand-Tsetlin procedure for the construction of orthogonal bases in Hermitean Clifford analysis},
In: T.E.\ Simos, G.\ Psihoyios, Ch.\ Tsitouras, {\em Numerical Analysis and Applied Mathematics}, AIP Conference Proceedings, Rhodes, Greece (2010).
%
\bibitem{GT2H} F.\ Brackx, H.\ De Schepper, R.\ L\'{a}vi\v{c}ka, V.\ Sou\v{c}ek, 
\textit{Orthogonal basis of Hermitean monogenic polynomials: an explicit construction in complex dimension $2$}, In: T.E.\ Simos, G.\ Psihoyios, Ch.\ Tsitouras (eds.), {\em Numerical Analysis and Applied Mathematics}, AIP Conference Proceedings, Rhodes, Greece (2010). 
%
 \bibitem{CacGueMal} I. Ca\c c\~ao, K. G\"urlebeck, H.R. Malonek, Special monogenic polynomials and $L_2$-approximation. Advances in Applied Clifford Algebras 2001; 11(S2):47--60.
%
\bibitem{CacGueBock} I. Ca\c c\~ao, K. G\"urlebeck, S. Bock. Complete orthonormal systems of spherical monogenics - a constructive approach. In Methods of Complex and
Clifford Analysis, Son LH, Tutschke W, Jain S (eds). Proceedings of ICAM, Hanoi, SAS International Publications, 2004.
%
 \bibitem{BockCacGue} I. Ca\c c\~ao, K. G\"urlebeck, S. Bock, On derivatives of spherical monogenics. Complex Variables and Elliptic Equations 2006; 51(811):847--869.
%
 \bibitem{CM06} I. Ca\c c\~ao and H.~R. Malonek, Remarks on some properties of monogenic polynomials,
ICNAAM 2006. International conference on numerical analysis and applied mathematics 2006 (T.E.
Simos, G. Psihoyios, and Ch. Tsitouras, eds.), Wiley-VCH, Weinheim, 2006, pp. 596-599.
%
 \bibitem{CM08} I. Ca\c c\~ao and H.~R. Malonek, On a complete set of hypercomplex Appell polynomials,
Proc. ICNAAM 2008, (T. E. Timos, G. Psihoyios, Ch. Tsitouras, Eds.), AIP Conference Proceedings
1048, 647-650.
%
 \bibitem{cac} I. Ca\c c\~ao, \textit{Constructive approximation by monogenic polynomials}, PhD thesis, Univ. Aveiro, 2004.
%
\bibitem{cno} J. Cnops, \textit{Reproducing kernels of spaces of vector valued monogenics}, Adv. appl. Clifford alg. \textbf{6} 1996 (2), 219-232.
%
\bibitem{R1} R. Delanghe, \textit{On homogeneous polynomial solutions of the Riesz system and their harmonic potentials},
Complex Var. Elliptic Equ.  \textbf{52}  (2007),  no. 10-11, 1047--1061.
%
\bibitem{R2} R. Delanghe, \textit{On homogeneous polynomial solutions of generalized Moisil-Th\'eodoresco systems
in Euclidean space}, CUBO \textbf{12} (2010), 145-167.
%
%
\bibitem{DSS} R. Delanghe, F. Sommen, V. Sou\v cek, \textit{Clifford Algebra and Spinor-valued Functions}, Mathematics and Its Applications 53, Kluwer Academic Publishers, 1992. 
%
\bibitem{DLS} R. Delanghe, R. L\'avi\v cka and V. Sou\v cek, \textit{On polynomial solutions of generalized Moisil-Th\'eodoresco systems and Hodge-de Rham systems}, to appear in Adv. appl. Clifford alg.\ (arXiv:0908.0842 [math.CV], 2009).
%
\bibitem{DLS2} R. Delanghe, R. L\'avi\v cka and V. Sou\v cek, \textit{The Fischer decomposition for Hodge-de Rham systems in Euclidean spaces}, preprint.
%
\bibitem{DLS3} R. Delanghe, R. L\'avi\v cka and V. Sou\v cek, \textit{The Howe duality for Hodge systems}, In: Proceedings of 18th International Conference on the Application of Computer Science and Mathematics in Architecture and Civil Engineering  (ed. K. Gürlebeck and C. Könke),
Bauhaus-Universität Weimar, Weimar, 2009.
%
\bibitem{FCM} M.~I. Falc\~ao, J.~F. Cruz and H.~R. Malonek, Remarks on the generation of monogenic
functions, Proc. of the 17-th International Conference on the Application of Computer Science and
Mathematics in Architecture and Civil Engineering, ISSN 1611-4086 (K. G\"urlebeck and C. K\"onke,
eds.), Bauhaus-University Weimar, 2006.
%
\bibitem{FM} M.~I. Falc\~ao and H.~R. Malonek, Generalized exponentials through Appell sets in $\bR^{n+1}$ and
Bessel functions, Numerical Analysis and Applied Mathematics (T.E. Simos, G. Psihoyios, and Ch.
Tsitouras, eds.), AIP Conference Proceedings, vol. 936, American Institute of Physics, 2007, pp.
750-753 (ISBN: 978-0-7354-0447-2).

\bibitem{GT} I.~M.~Gelfand and M.~L.~Tsetlin, {\it Finite-dimensional representations of groups of orthogonal
matrices}, Dokl. Akad. Nauk SSSR \textbf{71} (1950), 1017--1020 (Russian). English transl.\ in: I. M.
Gelfand, Collected papers, Vol.\ II, Springer-Verlag, Berlin, 1988, 657--661.
%
%
\bibitem{GM} J. E. Gilbert and M. A. M. Murray, \textit{Clifford Algebras and Dirac Operators in Harmonic Analysis},
Cambridge University Press, Cambridge, 1991.
%
\bibitem{GM10b}
K. G\"urlebeck and J. Morais, \textit{Real-Part Estimates for Solutions of the Riesz System in $\bR^3$}, 
to appear in Complex Variables.
%
\bibitem{Gurlebeck1999} K. G\"{u}rlebeck and H.~R. Malonek, A hypercomplex derivative of monogenic functions
in $\mathbb{R}^{n+1}$ and its applications, Complex Variables 39 (1999), 199--228.
%
\bibitem{NGue2009} N. G\"urlebeck, On Appell Sets and the Fueter-Sce Mapping, Advances in Applied Clifford Algebras 19 (2009), 51-61.
%
%
%
\bibitem{hom} Y. Homma, \textit{Spinor-valued and Clifford algebra-valued harmonic polynomials}, J. Geom. Phys. \textbf{37} (2001), 201-215.
%
\bibitem{lav_isaac09} R. L\'avi\v cka, \textit{The Fischer Decomposition for the $H$-action and Its Applications}, to appear.
%
\bibitem{lavSL2} R.~L\'avi\v cka, Canonical bases for sl(2,C)-modules of spherical monogenics in dimension 3, arXiv:1003.5587v2 [math.CV], 2010,
to appear in Arch. Math.
%
\bibitem{lavSL3} R.~L\'avi\v cka, 
On the Structure of Monogenic Multi-Vector Valued Polynomials, In: ICNAAM 2009, Rethymno, Crete, Greece, 18-22 September 2009 (eds. T. E. simos, G. Psihoyios and Ch. Tsitouras), AIP Conf. Proc. 1168 (2009)(793), pp. 793-796.
 %
 \bibitem{step2} R.~L\'avi\v cka, V.~Sou\v cek, P.~Van Lancker, Spherical monogenics: step two branching, preprint.
%

\bibitem{leu} H. Leutwiler, \textit{Quaternionic analysis in $\bR^3$ versus its hyperbolic modification}, In: F. Brackx, J. S. R. Chisholm and V. Sou\v cek (Eds.), Clifford analysis and its applications, Kluwer (2001), pp. 193-211.
%
\bibitem{CM07} H.~R. Malonek, M.~I. Falc\~ao, Special monogenic polynomials\,-\,properties and applications. In Numerical Analysis and Applied Mathematics (T.E.
Simos, G. Psihoyios, and Ch. Tsitouras, eds.), AIP Conference Proceedings, vol. 936. American Institute of
Physics: Melville, NY, 2007; 764--767.
%
\bibitem{Malonek1987} H.~R. Malonek, {Z}um {H}olomorphiebegriff in h\"{o}heren {D}imensionen, Habilitationsschrift. P\"{a}dagogische Hochschule Halle,
 1987.
%
 \bibitem{MS} I.~M. Mitelman and M.~V. Shapiro, Differentiation of
the Martinelli--Bochner integrals and the notion of hyperderivability. {\em Math. Nachr.} 172:
(1995), 211--238.
%
\bibitem{mol} A.~I.~Molev, {\it Gelfand-Tsetlin bases for classical Lie algebras},
in "Handbook of Algebra", Vol. 4, (M. Hazewinkel, Ed.), Elsevier, 2006, 109-170.
%
\bibitem{mor09} J. Morais, Approximation by homogeneous polynomial solutions of the Riesz system in $\bR^3$, PhD thesis, Bauhaus-Univ., Weimar, 2009.
%
\bibitem{som} F. Sommen, Spingroups and spherical means III, Rend. Circ. Mat. Palermo
(2) Suppl. No 1 (1989), 295-323.
%
\bibitem{van} P. Van Lancker, Spherical Monogenics: An Algebraic Approach, Adv. appl. Clifford alg. 19 (2009), 467-496.
%
\bibitem{vil} N. Ja. Vilenkin, 
\textit{Special Functions and Theory of Group Representations}, Izdat. Nauka,
Moscow, 1965 (Russian). English transl.\ in: Transl. Math. Monographs, Vol. 22, Amer. Math. Soc, Providence, R. I.,
1968.
%
\bibitem{zei} P. Zeitlinger, \textit{Beitr\"age zur Clifford Analysis und deren Modifikation}, PhD thesis, Univ. Erlangen, 2005 (German).
%
\end{thebibliography}
\end{document}